\documentclass[11pt]{article}
\usepackage{amsmath,amssymb,latexsym,graphicx}

\newtheorem{theo}{Theorem}[section]
\newtheorem{propo}[theo]{Proposition}

\input amssym.def
\newsymbol\rtimes 226F
\newfont{\nset}{msbm10}
\newcommand{\ns}[1]{\mbox{\nset #1}}

\def\Z{\ns{Z}}

\def\Cay{\mathop{\rm Cay }\nolimits}

\def\diam{\mathop{\rm diam }\nolimits}

\def\mod{\mathop{\rm mod }\nolimits}
\def\rank{\mathop{\rm rank }\nolimits}

\def\vec0{\mbox{\boldmath $0$}}

\begin{document}

\title{Comments on ``Extremal Cayley digraphs
of finite Abelian groups" [Intercon. Networks 12 (2011), no. 1-2, 125--135]
}

\author{M.A. Fiol\\
\\ {\small Departament de Matem\`atica Aplicada IV}
\\ {\small Universitat Polit\`ecnica de Catalunya}
\\ {\small Jordi Girona 1-3 , M\`odul C3, Campus Nord }
\\ {\small 08034 Barcelona, Catalonia (Spain); email: {\tt fiol@mat.upc.es}
}   }

\maketitle
\begin{abstract}
We comment on the paper ``Extremal Cayley digraphs
of finite Abelian groups" [Intercon. Networks 12 (2011), no. 1-2, 125--135]. In particular, we give some counterexamples to the results presented there, and provide a correct result for degree two.
\end{abstract}

\section{Introduction}
For the description of the problem, its applications, used notation, and the theoretical background, see, e.g. \cite{WoCo85,BeCoHs95,HsJi94,FiYeAlVa87}.

For some given positive numbers, $d$ (diameter) and $k$ (degree), the authors of \cite{MaSchJi11} consider the following numbers:
\begin{itemize}
\item
Let $m_{\ast}(d,k)$ the largest positive integer $m$ (number of vertices) such that there exists an $m$-element finite Abelian group $\Gamma$ and a $k$-element generating subset $A\subset\Gamma$ such that $\diam(\Cay(\Gamma,A))\le d$.
\item
Let $m(d,k)$ the largest positive integer $m$  such that there exists a cyclic  group $\Z_m$ and a $k$-element generating subset $A\subset\Z_m$ such that $\diam(\Cay(\Z_m,A))\le d$.
\end{itemize}

Such authors claim that, for any integer $d\ge 2$, Jia and Hsu \cite{HsJi94} proved that
\begin{equation}
\label{old-result}
m(d,2)=\left\lfloor \frac{d(d+4)}{3}\right\rfloor +1,
\end{equation}
but this was proved around ten years before by the author et al. \cite{MoFiFa85,FiYeAlVa87}.
Thus, in \cite{MoFiFa85}, the following value can be found:
\begin{equation}
\label{our-result}
m(d,2)=\left\lceil \frac{(d+2)^2}{3}\right\rceil -1,
\end{equation}
which is readily seen to be equivalent to \eqref{old-result}.
More generally, in Table I of \cite{FiYeAlVa87} some other optimal values are shown (minimizing the diameter for some fixed number of vertices). A part of such a table is shown next with the corresponding generating sets $\{a,b\}$ of the cyclic groups. (The values in boldface correspond to the ones given by \eqref{old-result} or \eqref{our-result}.

\medskip\begin{center}
\begin{tabular}{@{\,}cccc@{\,}}
$m(d,2)$ & $d$ & $a$ & $b(\mod m)$\\
\hline
$3x^2$ & $3x-1$ & $1$ & $3x-1$ \\
$3x^2+x$ & $3x-1$ & $1$ & $3x$ \\
{\boldmath $3x^2+2x$} & {\boldmath$3x-1$} & {\boldmath $1$} & {\boldmath $-3x$} \\
$3x^2+2x+1$ & $3x$ & $1$ & $3x+1$ \\
$3x^2+3x+1$ & $3x$ & $1$ & $3x+2$ \\
{\boldmath $3x^2+4x+1$} & {\boldmath $3x$} & {\boldmath $1$} & {\boldmath $-3x-2$} \\
$3x^2+4x+2$ & $3x+1$ & $1$ & $3x+3$ \\
$3x^2+5x+2$ & $3x+1$ & $1$ & $3x+4$ \\
{\boldmath $3x^2+6x+2$} & {\boldmath $3x+1$} & {\boldmath $1$} & {\boldmath $-3x+4$} \\
$(=3(x+1)^2-1)$ & & & \\
\hline
\end{tabular}
\end{center}
\medskip

Also, as a main result, Mask, Schneider, and Jia \cite[Th. 1.1]{MaSchJi11} claimed that, for any $d$ and $k$,
\begin{equation}
\label{noresult}
m_{\ast}(d,k)=m(d,k).
\end{equation}
However, as shown by the counterexamples in the following section, such a result cannot be true even for degree $k=2$.
This is due to an error in the proof of such a theorem. Namely, the first $r$ equalities in \cite[Th. 1.1]{MaSchJi11} should be understood modulo $m_j$:
$$
x_j=\sum_{i=1}^k c_ia_{ij}\quad (\mod m_j)\qquad\mbox{for}\ j=1,2,\ldots,r.
$$
Thus, without this condition, the following equality in \cite{MaSchJi11}, which should be modulo $m_{r-1}'=m_{r-1}m_r$, does not necessarily holds.

\section{Some counterexamples and a result}
In \cite{FiYeAlVa87} it was shown that for degree $k=|A|=2$, the minimum diameter $d$
of an Abelian group $\Gamma$ with $m$ vertices is $d_{\min}=\lceil\sqrt{3m}\rceil-2$ (see \cite[Eq. (9)]{FiYeAlVa87}). That is,
\begin{equation}
\label{Abel-bound}
m_{\ast}(d,2)\le \left\lfloor\frac{(d+2)^2}{3}\right\rfloor.
\end{equation}
In fact the upper bound is attained when $\Gamma=\Z_{3x}\times \Z_x$, with $x\ge 1$, and $A=\{(1,0),(-1,1)\}$, leading to a ($2$-regular) Cayley digraph on $m=3x^2$ vertices and diameter $d=3x-2$.
However, it can be shown that, when $x>1$,  $\rank \Gamma=2$, so that $\Gamma$ is not cyclic.
In this case, the best result is obtained with the cyclic group $\Z_m$ with $m=\frac{1}{3}(d+2)^2-1$ and generating set $A=\{a,b\}$, as shown in the following table.

\medskip\noindent
\begin{tabular}{@{\,}ccccccc@{\,}}
$k$ & $x$ & $d=3x-2$ & $m_{\ast}(d,2)=3x^2$ &$A\subset \Z_{3x}\times \Z_x$ & $m(d,2)=3x^2-1$ & $A\subset \Z_m$ \\
\hline
2 & 2 & 4 & $12$ & $\{(1,0),(-1,1)\}$ & $11$ & $\{1,3\}$ \\
2& 3 & 7 & $27$ &  $\{(1,0),(-1,1)\}$ & $26$ & $\{1,8\}$ \\
2& 4 & 10 & $48$ & $\{(1,0),(-1,1)\}$ &$47$ & $\{1,11\}$ \\
2& 5 & 13 & $75$ & $\{(1,0),(-1,1)\}$ &$74$ & $\{1,14\}$ \\
2& 6 & 16 & $108$ & $\{(1,0),(-1,1)\}$ &$107$ & $\{1,17\}$\\
\hline
\end{tabular}

\medskip\noindent
For other values of $m(d,2)$, see \cite[Table II]{FiYeAlVa87} or the results in \cite{EsAgFi93,AgFi95}.
In fact, from the results of these papers, and comparing the values of $m(d,2)$ in \eqref{our-result} with the upper bound for $m_{\ast}(d,2)$ in \eqref{Abel-bound}, one gets the following result for the case of degree $k=2$:
\begin{propo}
For any diameter $d\ge 2$,
\begin{equation}
m_{\ast}(d,2)=\left\{
\begin{array}{ll}
m(d,2)+1, & \mbox{if $d\equiv 1$\ $(\mod 3)$}, \\
m(d,2), & \mbox{otherwise.}
\end{array}
\right.
\end{equation}
\end{propo}

In the case of the above digraph $\Cay(\Z_{3x}\times \Z_x), \{(1,0),(-1,1)\})$, it can be shown that the two unique vertices at maximum distance $d=3x-2$ from the origin are
$(2x,x-1)$ and $(x,x-1)$.

Similar counterexamples can be given to prove that the extremal Cayley digraphs with respect to their average distance are not necessarily attained for cyclic groups (\cite[Th. 3.1]{MaSchJi11}).

\subsection*{Acknowledgements} Research supported by the
{\em Ministerio de Ciencia e Innovaci\'on}, Spain, and the
{\em European Regional Development Fund} under project MTM2011-28800-C02-01,
and the {\em Catalan Research Council} under project 2009SGR1387. The author wish to thank Hebert P\'erez-Ros\'es for calling my attention to, and supplying me with a copy of, the paper \cite{MaSchJi11}.


\begin{thebibliography}{99}

\bibitem{AgFi95}
F. Aguil\'{o}, M.A. Fiol,
An efficient algorithm to find optimal double loop networks,
{\it Discrete Math.} {\bf 138} (1995) 15--29.

\bibitem{BeCoHs95}
J.-C. Bermond, F. Comellas, and D.F. Hsu, Distributed loop computer networks: a survey, {\em J. Parallel Distribut. Comput.} {\bf 24} (1995) 2--10.

\bibitem{EsAgFi93}
P. Esqu\'{e}, F. Aguil\'{o}, M.A. Fiol,
Double commutative-step digraphs with minimum diameters,
{\it Discrete Math.} {\bf 114} (1993) 147--157.

\bibitem{HsJi94}
D.F. Hsu and X.D. Jia, Extremal problems in the construction of distributed loop networks, {\em SIAM J. Discrete Math.} {\bf 7} (1994) 57--71.

\bibitem{MaSchJi11}
A.G. Mask, J. Schneider, X. Jia, Extremal Cayley digraphs of finite Abelian groups,
{\em J. Intercon. Networks} {\bf 12} (2011), no. 1-2, 125--135.

\bibitem{MoFiFa85}
P. Morillo, M.A. Fiol, and J. F\`abrega, The diameter of directed graphs associated to plane tessellations, {\em Ars Combin.} {\bf 20-A} (1985) 17--27.

\bibitem{FiYeAlVa87}
M.A. Fiol, J.L.A. Yebra, I. Alegre, M. Valero, A discrete
optimization problem in local networks and data alignment, \emph{IEEE Trans.
Comput.} {\bf C-36} (1987) 702--713.

\bibitem{WoCo85}
C.K. Wong, D. Coppersmith, A combinatorial problem related to multinode memory organitzations,
\emph{J. Assoc. Comput. Machin.} \textbf{21} (1974) 392--402.


\end{thebibliography}
\end{document}